\documentclass[13pt, twoside, a4paper]{article}
\usepackage{amsfonts, amssymb, amsthm}
\usepackage[latin1]{inputenc}
\usepackage[intlimits,fleqn]{amsmath}
\usepackage[english]{babel}
\usepackage{graphicx} 

\include{epsf}

\newcommand{\la}{\langle}
\newcommand{\ra}{\rangle}

\newcommand{\suli}[2]{\sum\limits_{#1}^{#2}}
\newcommand{\ili}[2]{\int\limits_{#1}^{#2}}

\def\deg{\textrm{deg}}

\def\eh{\frac{1}{2}}

\def\be{\begin{eqnarray*}}
\def\bel{\begin{eqnarray}}
\def\ee{\end{eqnarray*}}
\def\eel{\end{eqnarray}}
\def\a{\alpha}
\def\b{\beta}
\def\g{\gamma}
\def\l{\lambda}
\def\w{\omega}

\def\M{\mathcal{M}}
\def\R{\mathbb{R}}
\def\N{\mathbb{N}}

\def\E{\mathbb{E}}
\def\P{\mathbb{P}}

\def\Cm{\mathcal{C}}
\def\Cn{\mathcal{C}_o}

\def\F{\mathcal{F}}

\def\Z{\mathbb{Z}}

\def\wh{\widehat}

\def\e{\epsilon}
\def\d{\delta}

\def\id{\mathbb{I}}

\def\Sp{\textrm{Tr}}

\def\nn{\nonumber}

\def\eqiv{\Leftrightarrow}

\def\vec{\left(\begin{array}{cc} }
\def\cev{\end{array} \right)}

\theoremstyle{plain}
\newtheorem{theo}{Theorem}[section]
\newtheorem{lemma}[theo]{Lemma}

\newtheorem{corollary}[theo]{Corollary}

\theoremstyle{definition}

\begin{document}\large
\begin{titlepage}

{\center
{\LARGE

\quad {\bf Bounds for the annealed return probability on large finite
  percolation graphs}}

\vspace{0.5cm}

 Florian Sobieczky \footnote{Mathematisches Institut, University of Jena,
   supported by FWF (Austrian Science Fund), project P18703.}\\
\vspace{0.5cm}}

{\bf Abstract:} Bounds for the expected return probability of the
delayed random walk on finite clusters of an invariant percolation on
transitive unimodular graphs are derived. They are particularly suited
for the case of critical Bernoulli percolation and the associated
heavy-tailed cluster size distributions.  The upper bound relies on
the fact that cartesian products of finite graphs with cycles of a
certain minimal size are Hamiltonian. For critical Bernoulli bond
percolation on the homogeneous tree this bound is sharp. The
asymptotic type of the expected return probability for large times $t$
in this case is of order $t^{-3/4}$.\\




\end{titlepage}  

\newpage

\section{Introduction}

\subsection{Context and Results}

This paper is about the expected return probability of the delayed
random walk on the finite clusters of percolation graphs with
heavy-tailed cluster size distributions (such as critical Bernoulli
percolation).\\

The asymptotics of the integrated density of states (IDS) of the graph
Laplacian on percolation subgraphs of the Euclidean lattice has
recently been studied in the subcritical phase by Kirsch and M\"uller
[\ref{kirsch}], and the supercritical phase by M\"uller and Stollmann
[\ref{muesto}]. The question of the IDS' asymptotics in the critical
phase was left open. For the two-dimensional Euclidean lattice, we
present upper and lower polynomial bounds (Theorem \ref{theo:2}) for
general invariant percolation. More generally, we find polynomial
bounds for the expected return probability on finite critical
percolation clusters on any planar transitive unimodular graph
(Theorem \ref{theo:3}). The upper estimates also hold in the
non-planar case. For homogeneous trees, this bound proves to be sharp
if the asymptotic type of decay of the cluster size' probability
density function is known (Theorem \ref{theo:lower}). Furthermore,
improved bounds for the number of open clusters per vertex
[\ref{grim_num}] in terms of the expected return probability are found
(Theorem \ref{theo:num}).\\

The method from which these bounds are derived are comparison theorems
for random walks on finite graphs. For the upper bound, the main idea
is the comparison of {\em all} the eigenvalues of the transition
matrices. Taking into account the whole spectrum instead of only the
spectral gap leads to an additional polynomially decreasing prefactor
in front of the exponentially converging return probability. For the
{\em expected} return probability an additional integration over all
finite random clusters is involved. As in critical percolation, the
corresponding cluster size distribution is heavy-tailed, i.e. integral
moments do not exist [\ref{timar}]. The result is a polynomial decay
in time. For this decay the additional prefactor is an essential
improvement.\\

The comparison theorem is obtained from the property of cartesian
products of finite graphs with maximum vertex-degree $\d$ and cycles
$C$ of size equal to $\d$ to be Hamiltonian [\ref{bapi}].  
This cycle exists due to Hamiltonicity. In addition to
this fact, we will use that the return probability of a continuous
time random walk on a finite cartesian product graph factorises into
the return probabilities on its factors. Since the return
probabilities are known on the cycle, this gives a bound for the
return probability on the original graph.
For the lower bound, we resort to a result by Boshier [\ref{boshier}]
about the isoperimetric number of a finite graph (see [\ref{mohar}]):
This is an upper bound for the isoperimetric number of graphs with
bounded genus.  For planar graphs, this gives us a bound of the
spectral gap from above by Cheeger's inequality.\\

\subsection{Delayed Random walk on finite graphs}\label{sec:drw}

We now recall some standard facts from finite random walk theory.  We
write $\N_0$ for $\{0, 1, 2, 3, ....\}$, and $\R_+:=[0,\infty)$. Since
we will assume $|\Cn|<\infty$, we will reserve subscript `o' for
objects defined in connection with {\em finite} graphs.\\

Let $G_o=\la V_o, E_o \ra$ be a finite simple graph, i.e. the
vertex-set $V_o$ has {\bf finite} cardinality and there are no
multiple edges in $E_o$, nor are they directed or have coinciding
incident vertices (`loops'). Let $\d$ be the maximal occurring degree,
i.e.  $\d:=\max\{\;\textrm{deg}(v)\;|\; v\in V_o\}$, where
$\textrm{deg}(v)=|\{w\in V_0\;|\;\{v, w\} \in E_o\}|$.\\

We define the discrete-time {\bf delayed random walk (DRW)} on $G_o$
to be the nearest neighbour random walk [\ref{woess}] with state space
$V_o$, some initial distribution $\nu\in \mathcal{M}_{+,1}(V_o)$, and
transition probabilities $P_{vw}:=(P(G_o))_{vw}$ with $v,w \in V_o$, and
\be 
\left(P(G_o)\right)_{v w}\;\;=\;\;\left\{\begin{array}{cr}
1/\d& \{v, w\} \in E_o, \\ 1\;-\;\deg(v)/\d &
v=w,\\0&\textrm{ otherwise.} \end{array} \right.  \ee 

 Recall that the transition probabilities of $v$ to $w$
after $n$ steps is given by the element of the matrix-power
$(P^n)_{vw}$, for all $v,w \in V_o$.\\

The continuous-time version of the delayed random walk with
coordinate-map $X_t$ is defined as the Markov-process on the
right-continuous $V_o$-valued functions depending on $t\in\R_+$, with
some initial distribution $\nu\in\M_{+,1}(V_o)$, and transition
probabilities
\bel
\P[X_t=w|X_0=v]\;\;=\;\;\left(e^{-t(1-P)}\right)_{v w},
\hspace{1cm}v, w \in V_o.\label{eq:contprob}
\eel

We note that $\left(e^{-t(1-P)}\right)_{v w}=\suli{n=0}{\infty}
\left(P^n\right)_{v w} \frac{t^n}{n!}e^{-t}$, and that $\left( P^n
\right)_{v w}$ is also the \-pro-bability of $X_t$ to reach $w$ from
$v$ {\em conditioned} on the event of there having been exactly $n$
jumps up to time $t$. The number $e^{-t}t^n/n!$ is the probability of
that event, which is also characterised by $t\in [t_n, t_{n+1})$, where
$t_n$ is the sum of $n$ independent exponentially distributed random
variables (`waiting times') with parameter 1. So
$\left(e^{-t(1-P)}\right)_{v w}=\sum_{n=0}^\infty\P[X_t=w|X_0=v, \;t\in[t_n,
t_{n+1})\,] \P[t\in[t_n, t_{n+1})]$,  (see [\ref{norris}]).\\

Finally, we note that choosing the initial distribution
$\nu\in\M_{+,1}(V_o)$ to be the {\bf uniform distribution},
i.e. $X_0\sim $UNIF$(\Cn)$, and $\nu(\{v\}) = 1/|V_o|$ gives the
return probability as the value of a normalised trace
\bel
\P[X_t=X_0]
\;=\;\suli{v\in V_o}{}\left( e^{-t(\id-P)}\right)_{v v}
\frac{1}{|V_o|}\;=\;\frac{1}{|V_o|}\Sp[e^{-t(\id-P)}],\label{eq:trace}
\eel

as $\P[X_t=X_0]=\sum_{v\in V_o}{}\P[X_t=X_0|X_0=v]\P[X_0=v]$, and
$\P[X_0=v]=\frac{1}{|V_o|}$.

\subsection{Invariant percolation on unimodular graphs}\label{sec:inv}

We now define the setting for which the results of section $2.1$ will 
be applied (see section $2.2$).\\

Let $G=\la V, E\ra$ be an {\bf infinite} simple (see above) graph,
which has a transitive, unimodular subgroup $\Gamma$ of the
automorphism group Aut$(G)$. `Transitive' means vertex-transitive,
here, i.e. for all $v, w \in V$, there is an automorphism $\g\in
\Gamma$, s.t. $w=\gamma(v)$. `Unimodular' means that the left Haar
measure of $\Gamma$ is the same as the right Haar measure. We call such
a graph a {\bf unimodular graph}.\\

A well-known result for unimodular graphs is the so called {\bf
  mass-transport-principle} (see [\ref{lyons}],[\ref{blps2}]).  It
says that for all $\Gamma$-diagonally invariant functions
($f(\gamma(v), \gamma(w))=f(v,w)$ for all $\gamma\in \Gamma$) it holds
that
\be
\suli{w\in V}{} f(v, w)\;\;=\;\;\suli{w\in V}{} f(w, v).
\ee

Let now $(\Omega, \F, \mu)$ be the probability space with $\Omega=2^E$
the two-valued functions on the edges and $\F=\otimes_E \F_o$ the
\mbox{product} $\sigma-$ algebra with $\F_o=\{\emptyset, \{0\}, \{1\},
\{0, 1\}\}$. On $\F$, we consider a probability distribution
$\mu:\F\to[0,1]$ with the property of $\Gamma$-invariance:
\be 
\mu(A)\;\;=\;\;\mu(\gamma(A)),\hspace{1cm}\textrm{for all }A\in \F,
\g\in\Gamma.  
\ee

In this way, for any fixed $\w\in\Omega$, we obtain a random subgraph
$G'(\w)\le G=\la V, E\ra$, of the form $G'(\w)=\la V, E'(\w) \ra$,
where 
\be 
E'(\w) \;\;=\;\; \{\;e\in E\;|\;\w(e)=1\;\}\;\;=\;\;\w^{-1}(\{1\}).  
\ee 

A subgraph of $G$ in which only edges are removed is called a {\em
  partial graph} of G. Therefore, with every $\w\in\Omega$, we
associate the random partial graph \\[-0.3cm]
\be
G'(\w)=\la V, \w^{-1}(\{1\})\ra.\\[-0.3cm]
\ee

We call the pair $\la G, \mu\ra$ an {\bf invariant percolation $\mu$
  on a unimodular graph G}.\\

We will now fix an arbitrary vertex $o\in V$, the `root', and look
for fixed $\w\in\Omega$ at the connected component of the
graph $G'(\w)$ which contains $o$, and call it $\Cn(\w)$.  Since we
will assume $|\Cn|<\infty$, we will be interested in invariant
percolation measures $\mu$ with $\mu$-almost surely finite connected
components, i.e.  
\be 
\mu(\{\w\in\Omega\;|\; |\Cn(\w)|<\infty\})\;=\;1.  
\ee
\newpage

{\bf Examples:} a.) {\em Bernoulli Percolation on the Euclidean
  Lattice}: $G=\la \Z^d, N.N.\ra$ (`Nearest Neighbours'), and $\mu$ is
the product measure on $\Omega$: $\mu=\otimes_{e\in E}\pi_e$, where
$\pi_e:\F_o\to [0,1]$, and $p=\pi(w(e)=1)\in[0,1]$, for all $e\in
E$. It is well-known that for sufficiently small $p$, the connected
components are a.s. finite (`subcritical regime').  Also, in the
`supercritical regime' or the `critical regime', for which
$\mu(|\Cn|=\infty)\;>\;0$, we may condition on the event $A:=\{\w\in
\Omega\;|\; |\Cn|<\infty \}$. The conditional measure
$\mu(\cdot|A)=\mu(\cdot\cap A)/\mu(A)$ is also $\Gamma$-invariant. It
is a celebrated result that Bernoulli bond-percolation has almost
surely finite clusters in the case $d=2$.\\

b.) {\em Bernoulli Percolation on homogeneous trees}. The Bernoulli
percolation measure $\mu$ on a homogeneous tree of degree $\d$ is
invariant under the action of any transitive subgroup of its
automorphism group. It is well-known (see [\ref{grimmett}], Chap. 10.1),
that for critical percolation on the binary tree, we have that the
$\P_\mu[|\Cn|\ge m] \sim m^{-\eh}$.\\

Now, we define the delayed random walk on a random partial
graph: Given $\w\in\Omega$, consider the finite subgraph of $G'(\w)$
induced by $\Cn(\w)$, i.e. (using a standard notation) consider 
\be
G_o(\w)\;\;:=\;\;G'(\w)|\Cn(\w).
\ee
As discussed in Section \ref{sec:drw} this induces a random finite random
walk with random state space $\Cn(\w)$, random initial distribution
$\nu^{(\w)}\in\M_{+,1}(\Cn(\w))$, and corresponding random return
probabilities
\begin{eqnarray}
P^{(\w)}_{v w}:=\left(P(G_o(\w)\right)_{v w},\label{eq:Pr}
\end{eqnarray}

which form a $|\Cn(\w)|\times|\Cn(\w)|$ matrix, where $|\Cn(\w)|$ is
$\mu$-a.s. finite.\\

The random continuous-time random walk is formed analogously to the
procedure of section 1.2, with $G_o=G_o(\w)$. Choosing
$\nu^{(\w)}\in\mathcal{M}_{+, 1}(\Cn(\w))$as the initial distribution
of the process to be the {\bf uniform distribution} on $\Cn(\w)$, 
the random continuous-time return probabilities turn out to be
(compare with (\ref{eq:trace}))
\be
\P[X^{(\w)}_t=X^{(\w)}_0]\;\;=\;\; \frac{1}{|\Cn(\w)|}\Sp[e^{-t(\id - P(\w))}],
\ee

where $P(\w)=((P^{(\w)})_{vw})$ with $v, w\in \Cn(\w)$ is the
transition probability matrix (\ref{eq:Pr}) of the random
discrete-time random walk on $\Cn(\w)$.\\

We are interested in the asymptotic corrections of the expectation
value of the return probabilities
\be P_t(o)\;\;:=\;\;\E_{\mu}\left[\frac{1}{|\Cn|} \Sp[e^{-t(\id -
      P)}] \right] \ee

for large values of the time $t>0$ from its limiting value, which is
given by $\E_\mu[1/|\Cn|]$.

\newpage

\section{Results}

We first present our estimates for finite graphs in Section 2.1, and
apply them in Section 2.2 to bound the expected return
probability. Sections 2.3 and 2.4 contain the applications concerning
the integrated density of states and the expected number of open
clusters per vertex.

\subsection{Bounds of the Return Probability on finite graphs}

\begin{theo}\label{theo:upperbound}
  Let $G_o=\la V_o, E_o\ra$ be a simple, finite, connected graph
  with $N$ vertices and largest degree $\d$.  Let $X_t$ be the delayed
  random walk on $G_o$, and $\beta_2$ the second-largest eigenvalue of
  its transition kernel.  For $X_0\sim $ UNIF$(V_o)$, $k\in
  \{1,..., N-2\}$ and $t>0$
\be 
i.)\;\;\;\;\;\;\P[X_t
  \;=\;X_0]\;&\le&\; \frac{1}{N}\;+\;
2\cdot\frac{k}{N}\,e^{-t(1-\beta_2)}\;+\;
\sqrt{\frac{\pi}{32}}\frac{\d \sqrt{\d+2}}{\sqrt{t}}
\exp\left(-\frac{32tk^2}{(\d+2)\d^2 N^2}\right),\\ 
ii.)\;\;\;\;\,\P[X_t
  \;=\;X_0]\;\;&\le&\;\; \frac{1}{N}\;+\; 2\cdot\frac{k
}{N}\,e^{-t(1-\beta_2)}\;+\;\frac{\d^2(\d+2)}{16t} \frac{N}{k}
\exp\left(-\frac{32 t k^2}{(\d+2)\d^2 N^2}\right).
\ee

iii.) If $G_o$ is also planar, and $N>288$, it holds for $t>0$
\be
 \;\;\;\;\;\;\P[X_t \;=\;X_0]\;\;&\ge&\;\;
\frac{1}{N}\;+\;\frac{\exp\left(-t K/\sqrt{N}\right)}{N},
\;\;\;\textrm{ with }K\;=\; 12\sqrt{2}\cdot\d.\\[-0.2cm]
\ee
\end{theo}

These bounds allow choosing an optimal value of $k$ if something
about the relation between $\beta_2$ and $N$ is known. If $k$ in
Theorem \ref{theo:upperbound}, i.) and ii.) is of the order of $N$,
the bound is qualitatively the same as the obvious estimate resulting
from using the Poincar\'e inequality $1-\beta_j\ge \d/(4 N^2)$ for 
all $j\in \{2, ..., N\}$ (see [\ref{laurent}], Chap. 3.2).
\\[-0.5cm]

\subsection{Annealed Return Probability on finite Percolation Subgraphs}

Let $P_t\,:=P_t(o)=\,\E_{\mu} \P[X_t = o \;|\; X_0 = o]\;$ denote the
expected return probability to the vertex $o$ of the continuous-time
delayed random walk on $\Cn(\w)$ at time $t\ge 0$.

\begin{theo}\label{theo:3}
For $\mu$ being any invariant percolation on a unimodular transitive
graph $G=\la V, E\ra$, let $A, B, a, b> 0$, with $b\le 2$ such that for
all $m\in \N$
\bel 
A \,m^{-a} \;\;\le \;\;\P_\mu[|\Cn|\ge m]\;\;\le\;\;
B\,m^{-b}.  \label{eq:assump1}
\eel
i.) Then with $C=5\left(4b/\d\right)^b(2\cdot 4^b \,+\,
\d(\d+2)/2)$, for all $\a$ with $0<\a<b$ and $t>0$
\be 
& P_t\;\;-\;\;\E_\mu\left[\frac{1}{|\Cn|}\right]
\;\;\le\;\;C\cdot \E_\mu[|\Cn|^{\a}]\;t^{-\eh(1+\a)}\;.&
\ee 
ii.) If $G$ is also assumed planar, and $K$ as in
Theorem \ref{theo:upperbound}, then for $t>\sqrt{288}$
\be 
D\cdot t^{-2a(1+1/b)} \;\;\le \;\; &P_t \;\;-\;\;
\E_\mu\left[\frac{1}{|\Cn|}\right]\;,&\hspace{0.6cm}\textrm{ where   
}\;\; D\;=\;e^{-K}\frac{A/2}{1+(2B/A)^{1/b}}\;.\\
\ee
\end{theo}

{\bf Remarks:} The folklore rule about easily obtained lower bounds
doesn't apply in this general setting of transitive graphs. The
quality of the argument of comparing the graph with the `host graph' $G$
on which the percolation is defined (see e.g. Lemma 2.2 in
[\ref{mathieu}]) generally gives poor results.  If for example $\Cn$
is the finite connected component containing the root of Bernoulli
percolation on a homogeneous tree with vertex-degree $\d$, then the
subtree of the homogeneous tree induced by a ball with radius equal to
that of $\Cn$ has typically a much {\em smaller}\, spectral gap.
Thus, it cannot be used for {\em lower} bounds of the return
probability. In the case of amenable graphs, however, this comparison
technique is successful (see e.g. [\ref{Tonci}] for results beyond the
Euclidean lattice). Furthermore, as upper bounds on the volume-growth
give lower bounds on the return probability (see e.g. [\ref{woess}],
Chap.  14.C), the lack of such a bound on the volume-growth under the
present assumptions comes at the cost of weaker results in Theorem
\ref{theo:3}, ii.). \\[-0.2cm]

Nevertheless, from the following discussion it will be seen that it is
for tree-like graphs $G$, for which the upper bounds perform well. The
upper bounds turn out better if few manipulations of the finite graph
in form of removals and additions of edges have to be undertaken to
retrieve a spanning cycle (we say that the graph is {\em similar} to
the spanning cycle). The proof (see Section \ref{sec:proofs}) involves
the comparison of the graph with that of a cycle having length
comparable to the graph's order (number of vertices). An example of
graphs for which this property may be likely to prevail is given by
finite subgraphs of the {\em incipient infinite cluster} of Bernoulli
percolation (see [\ref{grimmett}], Chap. 9.4). It occurs at the
critical retention probability $p_c$ under the additional condition of
being infinite [\ref{kesten}]. It is therefore of interest to compare
the expected return probability of the delayed random walk on the
incipient infinite cluster with the corresponding quantity on the
ordinary connected components of critical percolation to which Theorem
2.2 can be applied, as long as it has clusters at criticality which
are almost surely finite ([\ref{kesten3}], Theorem 2; [\ref{blps2}]):

\begin{corollary}\label{cor:z2_trees}
Consider $P_t$, the expected return probability of the delayed random
walk on finite percolation clusters of critical Bernoulli bond
percolation:

\begin{itemize}
\item[i.)]{For the $2$-dimensional Euclidean lattice, with
  $\a\in(0, 1/5]$ such that $\E_\mu[|\Cn|^\a]<\infty$, there is $C_2>0$
  such that for $\;t\ge 1$
\be
C_2^{-1} t^{-(1+\a^{-1})} \;\; \;\;\le \;\;\;\; P_t
\;\;-\;\;\E_\mu[1/|\Cn|]\;\;\;\;\le\;\; C_2 \,\E_\mu[|\Cn|^{\a}]\,
t^{-\eh(1+\a)}.
\ee
}
\item[ii.)]{For the homogeneous tree of degree $\d$, there is $\e>0$,
  and a constant $C_\d$ depending on $\e$, such that for $t\ge 1$
\bel
C_\d^{-1} t^{-3} \;\; \;\;\le \;\;\;\; P_t
\;\;-\;\;\E_\mu[1/|\Cn|]\;\;\;\;\le\;\; C_\d(\e) \,\E_\mu[|\Cn|^{\eh-2\e}]\,
t^{-\frac{3}{4}+\e}.\label{eq:trees}
\eel
}
\end{itemize}
\end{corollary}

{\bf Remarks:} It is easy to show that given $b>0$, the condition
$\P_\mu[|\Cn|\ge m]\;\;\le\;\; B\,m^{-b}$ for some $B>0$ implies
$\E_\mu[|\Cn|^\a]<\infty$ for all $\a\in\R$ such that $0<\a<b$. \\

The upper bound for the range of $\a$ in Corollary \ref{cor:z2_trees},
ii.) is a result by Kesten [\ref{kesten2}] (see also [\ref{grimmett}],
Table 10.1). The results obtained in Theorem 2.2 are valid for the
very general setting of any invariant percolation on a unimodular
transitive graph $G$, and therefore their quality varies strongly
depending on the structure of $G$ and the type of percolation measure
$\mu$. The $\a\le 1/5$ condition implies that the upper bound
Corollary \ref{cor:z2_trees} i.) for $P_t$ isn't stronger than $\sim
t^{-2/3}$, which would distinguish DRW on the finite critical
percolation cluster from the incipient infinite cluster if the
Alexander-Orbach conjecture would be true.  However, it isn't believed
that the this conjecture holds for the Euclidean lattice in dimensions
$d\le 6$ ([\ref{hughes}], Chap. 7.4.4).\\[-0.1cm]

The situation with Corollary \ref{cor:z2_trees}, ii.) is different. It
is clear from Lemma 1.6 of [\ref{sob}], that DRW and the simple random
walk SRW on any finite subgraph of an infinite graph of polynomial
growth have the same decay-type of the expected or quenched return
probability, as long as the maximum vertex-degree is uniformly
bounded.  Kozma and Nachmias (see Theorem 1.2 and 1.3 in [\ref{kona}])
have shown that the volume growth of the incipient infinite cluster in
high dimensional Euclidean lattices is almost surely polynomial. The
same follows from Lemma 2.2 of Barlow and Kumagai [\ref{baku}] for
homogeneous trees. Both of these cases are percolation models on
transitive graphs with uniformly bounded vertex-degree. Corollary
\ref{cor:z2_trees} is therefore interesting when compared with the
results obtained in [\ref{kona}] and [\ref{baku}] for the asymptotics
of the simple random walk on the incipient infinite cluster on
trees. It is proved there that the expected return probability is -
regardless of the degree $\d$ - of the order of $t^{-2/3}$. (That the
so called spectral dimension $-2 \lim_n \log \P_o[X_n=o]/\log n$ is
equal to $-4/3$ is known as the Alexander-Orbach conjecture
[\ref{alob}], proven for homogeneous trees [\ref{baku}], and Euclidean
lattices for high dimensions [\ref{kona}].)  Since (\ref{eq:trees})
represents an upper bound for $P_t-\E_\mu[1/|\Cn|]$ that can be chosen
to have an exponent arbitrarily close to $-3/4$, it proves that the
expected return probability at criticality on ordinary finite
percolation clusters displays a different asymptotic decay towards
its limit  than on the incipient infinite
cluster.\\[-0.3cm]

We expect the upper bound Theorem 2.2, i.) to be a good approximation
when $G$ is similar to a homogeneous tree and $\P_\mu[|\Cn|=m]$ is
polynomially decreasing in $m$:
\begin{theo}\label{theo:lower} 
Let $G$ be the homogeneous tree of degree $\d$ and $\mu$ an invariant
percolation on $G$ obeying assumption (\ref{eq:assump1}), and $A \le
\P_\mu[|\Cn|=m]m^{a+1}$ for all $m\in \N$. Then there is $c>0$, such
that for all $t>0$
\bel
P_t(o)-\E_\mu[1/|\Cn|]\;\;\ge\;\; c\,t^{-\eh(1+a)}.\label{eq:lower}
\eel
\end{theo}

We conclude the discussion of our results by the following tight
estimate for independent percolation on the homogeneous tree:

\begin{corollary}\label{cor:three_fourths} For critical Bernoulli 
bond percolation on the homogeneous tree
\be
\lim\limits_{t\to\infty} \frac{\log(\;P_t(o)\;-\;\E_\mu[1/|\Cn|]\;)}{\log
  t} \;\;=\;\;-\frac{3}{4}.
\ee
\end{corollary}

These findings allow to conclude that the observation by Kirsch and
M\"uller [\ref{kirsch}] of the predominance of path-like clusters also
determines the asymptotics of critical percolation in the present
case. The difference of (\ref{eq:lower}) over subcritical percolation
considered in [\ref{kirsch}] consists of the necessity to include, in
addition to the `linear' clusters [see Remark 1.15, iii.)  in
  [\ref{kirsch}]), the larger class of clusters $\Cn$ which have
  diameters $D$ comparable to the cluster's size $|\Cn|$.\\

The fact that path-like clusters are the dominating structures for the
large-time asymptotics in the case of trees is also illustrated by the
following `heuristic', but wrong argument: Suppose that for a given
realisation $\w\in\Omega$ at time $t>0$, the cluster size $|\Cn|$ is
larger than $t^{2/3}$. One might guess that up to this time the Markov
chain hasn't equilibrated and this cluster contributes significantly
in the averaging over $\P[X_t=X_0]- 1/|\Cn|$.  For times larger than
$t^{2/3}$ one then assumes that $\P[X_t=X_0]\sim 1/|\Cn|$. Assuming
further that up to the time of equilibriation the return probability
on the finite cluster typically decays just like on the incipient
infinite cluster, namely like $t^{-2/3}$ (see [\ref{baku}], Theorem
1.4), then by using $\P_\mu[|\Cn| = m]\sim m^{-3/2}$, one arrives at
the following rough estimate: 
\be
P_t(o)-\E_\mu\left[\frac{1}{|\Cn|}\right] \;\sim\; \sum\limits_{m \ge
  t^{2/3}}\left(t^{-2/3}-\frac{1}{m}\right) \,m^{-3/2} \;\sim\; t^{-1}, 
\ee

where the first $\sim$ (meaning `of this order, for large t') follows
from assuming that only unsignificant terms are neglected. This, however,
contradicts Corollary \ref{cor:three_fourths}.\\

The reason for restricting the considered clusters to sizes of at
least $t^{2/3}$ in this argument comes from the idea that because the
characteristic asymptotic decay of the random walk on the incipient
infinite cluster is $t^{-2/3}$ the random walk on smaller clusters
will have already reached equilibrium, and all of the significant
contributions to $(t^{-2/3}-1/|\Cn|)_+$ are accounted for. It is
therein implicitly assumed, that the typical decay on clusters of
smaller size before equilibriation is also $\sim t^{-2/3}$. However,
the path-like clusters have a characteristic heat-kernel decay towards
$1/|\Cn|$ of order $t^{-1/2}$ instead of $t^{-2/3}$ (see Part ii. and
iii. in the proof of Theorem \ref{theo:lower}). And so our result
shows that the regime of cluster sizes between $t^{1/2}$ and $t^{3/2}$
plays the dominant part in the averaging for Bernoulli percolation on
the homogeneous tree.\\

The reason why these contributions are not relevant in the case of the
infinite incipient cluster follows from the results of Barlow and
Kumagai [\ref{baku}]: According to their Lemmata 2.2 and 2.3 the
incipient infinite cluster on homogeneous trees has realisations
which, if restricted to subtrees with radius $n$, typically have a
size of order $n^2$ , so that the diameter ($\sim n$) is never a
positive fraction of the cluster size. From this it becomes apparent
that the main characteristic responsible for the stronger decay of the
upper bounds in Theorem \ref{theo:3} is the existence of a significant
fraction of finite clusters with diameter comparable to their size.

\subsection{Integrated density of states for $\Z^2$}

Let $\mu$ be an invariant bond percolation on the $2$-dimensional
Euclidean lattice $G=\la \Z^2, N.N.\ra$ with a $\mu$-a.s. finite
percolation cluster $\Cn$ having a size distribution obeying
(\ref{eq:assump1}). Let $\a\in(0, 2)$ such that
$\E_\mu[|\Cn|^\a]<\infty$.\\

Let $N(E)$ be the {\bf integrated density of states} (IDS) of the
graph Laplacian $L(\w)$ belonging to the percolation subgraphs
$G'(\w)$. This means for $\Lambda_N=\{-N+1, ..., N\}^2\subset\Z^2$ the
limit
\be 
N(E)\;\;=\;\;\lim\limits_{N\to\infty} \,\frac{1}{|\Lambda_N|}\#\{
\lambda \textrm{ eigenvalue of } L_{\Lambda_N}(\w) \le E\}
\ee

exists, where $L_{\Lambda_N}(\w)$ is the graph Laplacian of the finite
induced subgraph $G'(\w)|\Lambda_N$ (see e.g. [\ref{kirsch}], Lemma
1.12). 

\begin{theo}\label{theo:2}
There is $C_3 > 0$, s.t. for $E>0$ sufficiently
small the integrated density of states $E\mapsto N(E)$ of the graph
Laplacian obeys 
\be
 C_3^{-1} \frac{E^{1+1/\a}}{(\log1/E)^{1+1/\a}}
\;\;\;\;\le\;\;\;\; N_N(E)\;\;-\;\;N_N(0) \;\;\;\;\le\;\;\;\; C_3 \;
E^{\eh(1+\a)}.\\[-0.3cm] 
\ee
\end{theo}

{\bf Remarks:} This shows that independently of the vertex-degree $\d$,
the type of the asymptotics of the integrated density of states for
small values of $E>0$ is polynomial and only depends on the decay of
the cluster size distribution. By comparison with Theorem 1.14 of
[\ref{kirsch}], by which for subcritical percolation on the Euclidean
lattice in any dimension ($d=\d/2$)
\be
\exp(-\a_-/\sqrt{E})\;\;\le\;\;N(E)\;-\;N(0)\;\;\le\;\;
\exp(-\a_+/\sqrt{E}), 
\ee

for some $\a_-, \a_+ > 0$, and $E>0$ sufficiently close to zero, it is
seen that observation of the type of asymptotics of the IDS for small
energies suffices to decide about whether the finite random cluster of
the origin is generated with a critical, or subcritical percolation
measure.\\[-0.7cm]
 
\subsection{Number of open clusters per vertex}

A central theme in percolation theory on the Euclidean lattice $G=\la
\Z^d, N.N.\ra$ is the so called {\bf number of open clusters
  per vertex}. Given a finite box $\Lambda_N=\{ -N+1, ..., N\}^d$, and
the number $M_N(\w)$ of connected components of the induced subgraph
$G'(\w)|\Lambda_N$, the $\mu$-a.s. existence of the limit
\be
\kappa(p)=\lim\limits_{N\to\infty} \frac{M_N(\w)}{|\Lambda_N|}
\ee

and its almost sure independence of $\w\in\Omega$ has been shown by
Grimmett [\ref{grim_num}]. Its value equals
$\kappa(p)\;\;=\;\;\E_\mu[1/|\Cn|]$ (see
[\ref{grimmett}], (4.18) ).  Note that the number $1/\Cn(\w)$ is the
value of the density of the uniform distribution on $\Cn(\w)$. \\

In [\ref{grim_num}] there are upper and lower bounds for $\kappa(p)$
(there it is defined by $\E_\mu[1/|\Cn|] - \mu[|\Cn|=1]$) in the case
of Bernoulli percolation on the Euclidean lattice.  They entail
expansions which are converging slowly in the regime of the retention
probability $p$ being close to the critical value.  We present the
consequences of our bounds in terms of the expected cluster size
$\chi(p)=\E_\mu[|\Cn|]$:

\begin{theo}\label{theo:num}
Let $\mu$ be subcritical Bernoulli bond percolation on the 
$d$-dimensional Euclidean lattice $G=\la \Z^d, N.N.\ra$ with 
almost surely finite connected components. Let
  $\chi(p)=\E_\mu[|\Cn|]$. Then, for \;$t>0$
\bel
P_t\;\;-\;\;c\, \frac{\chi(p)}{t}\;\;\;\le\;\;\;\kappa(p)
\;\;\;\le\;\;\;P_t\,.
\label{eq:grim_num}
\eel
\;\;\;\;with
$c=\min\{\;\eh(d^3 + d^2 + 4),\;\; \frac{20}{d}(4+d(d+1))\}$.\\
\end{theo}

{\bf Remarks:} The power of the method for the upper bound (mainly due
to Lemma \ref{lemma:cart}) becomes visible if one compares Theorem
\ref{theo:num} with the simple bound obtained by using Poincar\'e's
inequality for $\lambda$, together with $\lambda \le 1-\beta_j$, for
$j\ge 2$: In this case $P_t\;-\;\kappa(p)\;\;\le\;\;
\E_{\mu}[e^{-t\d/(4|\Cn|^2)}]$ instead of (\ref{eq:grim_num}) which yields
for \;$t>0$
\be 
P_t\;\;-\;\;\frac{2}{d}\cdot\frac{\E_{\mu}[|\Cn|^2]}{ t}
\;\;\;\le\;\;\;\kappa(p)\;\;\;\le\;\;\;P_t.
\ee

The constant in front of the term $t^{-1}$ includes the second moment
of the cluster size, while in (\ref{eq:grim_num}) only the first
moment appears.

\section{Proofs}
\subsection{Auxiliary results}

The proofs rest on the theory of infinite unimodular transitive graphs
[\ref{blps2}]. `Unimodularity' of a graph refers to the existence of a
vertex-transitive subgroup of the automorphism group of the graph.

\begin{lemma}\label{lemma:mtp} Let $G$ be an infinite unimodular vertex-transitive
graph, an $\mu$ an invariant percolation measure on $G$.  If
\;$\E_{\mu}$ refers to the integration of the expected value over all
partial graphs $\w \in \Omega$,
\bel 
\E_{\mu}\left[\P_o[X_t
    \;=\;o]\right]\;\;\;=\;\;\;\E_{\mu}\left[\P[X_t \;=\;X_0]\right].
\label{eq:uni}
\eel
\end{lemma}

{\em Proof:} (see [\ref{sob}] for a detailed discussion) Let $\Cm_v$
be the connected component of $H(\w)$ containing the vertex $v\in
V$. Since the Euclidean lattice is a graph with a unimodular group of
automorphisms, by the mass-transport-principle [\ref{blps2},
  \ref{lyons}], the left-hand side of (\ref{eq:uni}) equals
\be 
\suli{v\in V}{}\E_{\mu} \left[ \P_o[X_t \;=\;o] \frac{\chi_{{}_{\{v
        \in\Cn\}}}}{|\Cn|}\right] = \suli{v\in V}{}\E_{\mu} \left[ \P_v[X_t
  \;=\;v] \frac{\chi_{{}_{\{o \in\Cm_{v}\}}}}{|\Cm_{v}|}\right] = \suli{v\in
  V}{}\E_{\mu} \left[ \P_v[X_t \;=\;v] \frac{\chi_{{}_{\{v
        \in\Cn\}}}}{|\Cn|}\right]
\ee 

since $v\in \Cn \;\eqiv\; o\in \Cm_v$, which equals the
right-hand side of(\ref{eq:uni}).  \hfill\qed\\

\begin{lemma}\label{lemma:first} For $N>3$ and $k\in\{1, ..., N-2\}$,  
let $I_t(k,N):= \suli{j=k+1}{N-1} e^{-t(1-\cos \pi
  \frac{j}{N})}$. Then \bel i.)  \;\;\;
\;\;\;I_t(k,N)\;\;\;&\le&\;\;\;
\eh\sqrt{\frac{\pi}{2}}\frac{N}{\sqrt{t}}\,
e^{-{2t}\frac{k^2}{N^2}},\label{eq:teclem1}\;\;\;\;\textrm{ and
}\\ ii.)\;\; \;\;\;\,I_t(k,N) \;\;\; &\le& \;\;\; \eh\frac{N^2}{k
  t}\,e^{-2t\frac{k^2}{N^2}}.
\label{eq:teclem2}
\eel
\end{lemma}

{\em Proof:} From $\cos \pi x \le 1-2 x^2$, if $x\in [0,1]$, we obtain by
following [\ref{laurent}], (Ex. 2.1.1)
\bel \suli{j=k+1}{N-1}e^{-t(1-\cos\pi
  \frac{j}{N})}\;&\le&\;\ili{k}{\infty} e^{-2t\frac{x^2}{N^2}}dx\;
=\; \frac{N}{\sqrt{2t}}\ili{\sqrt{2t}k/N}{\infty}e^{-y^2}dy
\label{eq:teclem3} \\&\le&\; \frac{N}{\sqrt{2t}} e^{-{2t}\frac{k^2}{N^2}}
\ili{0}{\infty} e^{-y^2-2\sqrt{2t}\frac{k}{N}y}dy\; \le\;  
\frac{N}{\sqrt{2t}} e^{-{2t}\frac{k^2}{N^2}} \ili{0}{\infty} e^{-y^2}dy, \nn
\eel

which proves (\ref{eq:teclem1}). Moreover, we have

\be \ili{z}{\infty}e^{-u^2}du\;\;=\;\; \eh \ili{z^2}{\infty}
e^{-y}\frac{dy}{\sqrt{y}}\;\;=\;\;\eh\ili{0}{\infty}e^{-(y+z^2)}
\frac{dy}{\sqrt{y+z^2}} \;\;\le\;\; \frac{e^{-z^2}}{2z}\ili{0}{\infty}
e^{-y} dy.  \ee

Applying this inequality to the right-hand side of (\ref{eq:teclem3})
with $z=\frac{\sqrt{2t}k}{N}$ gives (\ref{eq:teclem2}).\hfill \qed\\

\begin{lemma}\label{lemma:cart}
  Let $\wh{G}=G_X\square G_Y$ be the cartesian product of the simple,
  connected, finite graphs $G_X, G_Y$. Let $\wh{X}_t$ be the
  continuous-time delayed random walk on $\wh{G}$ with uniform initial
  distribution on the vertices of $\wh{G}$.  Let $X_t$ and $Y_t$ be
  the continuous-time delayed random walk on $G_X$ and $G_Y$, also
  with uniform initial distribution on the vertex-sets of $G_X$ and
  $G_Y$, respectively . Then 
\be
  \P[\,\wh{X}_{2t}\;=\;\wh{X}_0\,]\;=\;\P[\,X_t\;=\;X_0\,]
\cdot\P[\,Y_t\;=\;Y_0\,].\\ 
\ee
\end{lemma}

{\em Proof:} Let $N=|V(G_X)|$, and $M=|V(G_Y)|$. Let
$P_X$ and $P_Y$ be the transition kernels of $X_t$ and $Y_t$,
respectively.  For the delayed random walk on $\wh{G}$, with equal
transition weights across edges of type $\{ \la x, v\ra, \la y,
v\ra\}$, and $\{ \la x, v\ra, \la x, w\ra\}$ (where $x,y \in V(G_X)$,
and $v,w\in V(G_Y)$), the transition kernel is given by $\eh(P_X\otimes
\id\;+\;\id\otimes P_Y)$ (see [\ref{woess}], Chap. 18). Therefore, 
\be 
\P[\wh{X}_{2t}=\wh{X}_0]\;&=&\;\frac{1}{N\cdot M}
 \Sp[e^{-2t(\id-\eh(P_X\otimes\id \;+\;\id\otimes P_Y))}]=
\frac{1}{N\cdot M}\Sp[e^{-t(\id - P_X)}\otimes e^{-t(\id - P_Y)}] \\ &=&
  \frac{1}{N}\Sp[e^{-t(\id - P_X)}]\;\frac{1}{M}\Sp[e^{-t(\id -
      P_Y)}] = \P[X_{t}=X_0]\;\P[Y_{t}=Y_0].\;\;\;\qed 
\ee

{\bf Remark:} This auxiliary result can also be derived by using the
fact that the sum of two independent Poisson processes is also a
Poisson process, however with rate equal to the sum of the two
components' rates (see e.g. [\ref{norris}], Theorem 2.4.4).

\begin{lemma}\label{lemma:heavytail}
Let $\phi:\N_0 \to \R_+$, s.t. $\suli{k=0}{\infty} \phi(k)=1$ with
$\Phi(m):= \suli{k=m}{\infty} \phi(k)$.  Let there exist $A, B, a, b
\in \R_+$ such that
$
\frac{A}{m^a}   \;\;\le\;\; \Phi(m) \;\;\le\;\; \frac{B}{m^b}
$
for all $m \in \N$. Then
\be
\suli{k=m}{\infty} \frac{1}{k} \,\phi(k) \;\;\ge\;\;
\frac{C}{m^{a(1+1/b)}}, && \textrm{ with } \; \;\;C=\;
\frac{(A/2)^{1-1/b}}{B^{1/b}}.
\ee
\end{lemma}

{\em Proof:} \be \suli{k=m}{\infty} \frac{1}{k}\,\phi(k) \;\;\ge \;\;
\suli{k=m}{L} \frac{1}{k}\,\phi(k) \;\;\ge \;\;
\frac{1}{L}\left(\Phi(m)\;-\;\Phi(L+1)\right)) \;\;\ge \;\;
\frac{1}{L}\left( \frac{A}{m^a} \;\;-\;\;\frac{B}{(L+1)^b}\right).
\ee We set $\tilde{L}>0$ to be the real value $L$, such that the
parentheses on the right-hand side are exactly $\eh\cdot A/m^a$,
i.e. $\tilde{L}:=\left(\frac{2 B}{A} \right)^{1/b} m^{a/b}$. Now, by
defining $L_-:=\lfloor\tilde{L}\rfloor$ and $L_+:=L_-+1$, we have as a
lower bound for the right-hand side 
\be 
\frac{1}{L_-} \left( \frac{A}{m^a}  \;\;-\;\;\frac{B}{(L_+)^b}\right) 
\;\;\;\ge\;\;\; \frac{1}{\tilde{L}}\left( \frac{A}{m^a}
  \;\;-\;\;\frac{B}{{\tilde{L}}^{{}^{\textrm{\small {\em b}}}}}\right)
\;\;\;\ge\;\;\; \frac{1}{m^{a/b}\,\left(\frac{2B}{A} \right)^{1/b}} 
\cdot \frac{A}{2\,m^a}.\;\;\;\qed
\ee

\subsection{Proofs of main results}\label{sec:proofs}

{\em Theorem \ref{theo:upperbound}; Upper bounds}: By the
Theorem of [\ref{bapi}] (see also the discussion in [\ref{dipapo}])
the cartesian product $\wh{G}:=G_o\square C_{\d}$ is Hamiltonian. Let
$Y_t$ be the continuous-time delayed random walk on the cycle $C_\d$
of order $\d$, with transition kernel $P_Y$. Since
$1/\d\;\le\;\P[Y_t=Y_0]=(1/\d)\Sp \exp(-t(\id-P_Y))\;\le\;1$, and from
Lemma \ref{lemma:cart} it follows
\be 
\P[\wh{X}_{2t}\;=\; \wh{X}_0] \le \P[X_t \;=\;X_0]\;\le \; \d
\,\cdot\,\P[\wh{X}_{2t}\;=\;\wh{X}_0], 
\ee

where $\wh{X}_t$ is the continuous-time delayed random walk on $\wh{G}$. 
By Theorem 1 in [\ref{heuvel}], the eigenvalues of the transition
kernel $\wh{P}$ of $\wh{X}_t$ can be compared with the eigenvalues of
the delayed random walk on $C_{\d N}$; namely, 
\bel 
\wh{\b}_{j}\;\;\;\le \;\;\;1-\frac{2}{\d+2}\left(1-\cos 2\pi\,
\frac{j-1}{\d N}\right),\;\;\;\;\;\;\;\;\;(j\in \{1, ..., \d
N\}),\label{eq:cosinus} 
\eel

where $1=\wh{\beta}_{1}> \wh{\beta}_2 \ge \wh{\beta}_3 \ge
\wh{\beta}_4 \ge ...  \ge \wh{\beta}_{\d N}$, and $N=|V(G_o)|$. The
factor $2/(\d+2)$ in front of the parentheses results from the
regularisation with loops, characteristic of the delayed random walk
on a graph ($\wh{G}$) with maximal degree $\d+2$, where the extra $2$
comes from taking the cartesian product with $C_\d$ (see
[\ref{sob}]). Note, the eigenvalue of $\wh{P}$ can also be enumerated
differently: $\{\wh{\beta}_j\}_{j=1}^{\d N} = \{\wh{\beta}_{j,
  l}\}_{j,l=1}^{N,\d}$, where $\wh{\b}_{j, l}=
\eh(\beta_j\;+\;\cos(2\pi(l-1)/\d)$, with $j\in \{1, ..., N\}$ and
$l\in\{1, ..., \d\}$. From (\ref{eq:trace}), we have
$\P[\wh{X}_{2t}\;=\;\wh{X}_0] = \frac{1}{\d N}
\Sp[e^{-2t(1-\wh{P})}]$, so

\bel \d\, \P[\wh{X}_{2t}\;=\;\wh{X}_0] &=&
\frac{1}{N}\suli{j=1}{N}\suli{i=1}{\d} e^{-2t(1-\eh(\beta_j +
  \cos(2\pi(i-1)/\d)))}\nn\\ &\le& \frac{1}{N}(1 \;+\; 2\cdot k\,
e^{-t(1-\beta_2)}) \;+\; \frac{1}{N} \suli{j=k+2}{\d N-k}
e^{-2t(\frac{2}{\d+2}(1-\cos2\pi\frac{j-1}{\d N}))}\nn\\ &\le&
\frac{1}{N}(1\;+\;2\cdot k\, e^{-t(1-\beta_2)})+ \frac{2}{N} \suli{
  j=k+1}{\lfloor \d \frac{N}{2}\rfloor - 1}
e^{-\frac{4t}{\d+2}(1-\cos2\pi\frac{j}{\d N})}.\label{eq:twoterms}
\eel

The first inequality follows from bounding the first $2k$ eigenvalues
of $\wh{P}$ less than one from above by $\beta_2=\wh{\beta}_{2,1}$,
and from (\ref{eq:cosinus}), giving that the $n$-th largest
element of $\{ \wh{\beta}_{j, l}\}_{j,\,l =1}^{N, \d}$ is less than
the $n$-th largest eigenvalue of DRW on $C_{\d N}$, which however is
only applied to $n>2k+1$.  The second inequality follows from the
symmetry of the cosine-function, and an index-shift, with equality if
$\d N$ is even.  Since $I_t(\cdot,\cdot)$ is monotone in the second
argument, the claim follows from applying Lemma \ref{lemma:first} i.)
and ii.) to $I_t(4t/(\d+2),\d N/2)$.  \hfill\qed \\[-0.2cm]

{\bf Remark:} (Theorem \ref{theo:upperbound}) For $1-\beta_2$ we have
the standard lower bound given by the Poincar\'e inequality.  The
delayed random walk has the same spectrum as the simple random walk on
the path `decorated' with loops to yield a regular graph of degree
$\d$ [\ref{sob}]. In particular, $1-\beta_2\ge
1-(1-2/\d(1-\cos(\pi/N)))\ge 4/(\d N^2)$, by $\cos\pi x\le 1-2x^2$ for
$x\in[0,1]$. If $k\in \{1, ..., N-2\}$ in (\ref{eq:twoterms}) is
chosen such that $\frac{4}{\d N^2} \le \frac{32 k^2}{(\d+2)\d^2N^2}$,
or, equivalently $k^2\ge \d(\d+2)/8$, then the first exponential term
$\exp(-t(1-\beta_2))$ has weaker decay than the second. We see this is
the case for a number $k$ independent of $N$. Therefore, provided that
$N$ is sufficiently large, even if nothing else is known about
$\beta_2$, Theorem \ref{theo:upperbound} is an improvement over simply
using $\beta_2 \ge \beta_j$ for $j\ge 2$ and the Poincar\'e inequality
for $1-\beta_2$, which would be the bound corresponding to $k=N-1$ and
the second term in (\ref{eq:twoterms}) vanishing. \\

{\em Theorem \ref{theo:upperbound}; Lower bound}: 
For a given finite simple graph $G_o=\la V_o, E_o\ra$, let $I(G_o)$ be
the isoperimetric number (or `Cheeger-constant') of $G_o$, defined by
\be
I(G_o)\;\;\;=\;\;\;\min\limits_{{A\subset V_o}\;:\; {|A|\le \eh
    |V_o|}} \frac{|\partial_{G_o} A|}{|A|}\;, 
\ee 
where $\partial_{G_o}A=\{ \{k, l\} \in E_o\;|\; k\in A, l\notin A\}$
is the {\em edge-boundary} of $A$ in $G_o$, and $|A|=\#A$ denotes
cardinality of the finite set $A$.\\

By a theorem of A.G. Boshier [\ref{boshier}] the isoperimetric number
$I$ for graphs with genus bounded by $g$ obeys $I\le
3\,\d(g+2)/(\sqrt{|V_o|/2}-3(g+2))$ if $|V_o|>18(g+2)^2$ (see
[\ref{mohar}] for a discussion).  From this result it holds that if
$|V_o|\ge 4\cdot 72$ and $G_o$ a planar finite graphs (for which
$g=0$!) that
\bel 
I \;\;\;\le\;\;\; K/\sqrt{|V_o|},\;\;\hspace{1cm}\textrm{ with }
K\;=\;12\sqrt{2}\cdot\delta.  \label{eq:boshier}
\eel

By Cheeger's inequality (see [\ref{laurent}], Lemma 3.3.7), the
spectral gap $\lambda= \eh\min_{v\neq \textrm{const}} $ $(v,
(1-P) v)/(v,v)\;=\;1 -\beta_2$ for the delayed random walk with
transition probability matrix $P$ can be estimated from above,
\be 
\lambda \;\;\le\;\; I.
\ee

By (\ref{eq:boshier}) this implies a lower bound on the return
probability of the continuous-time delayed random walk for planar
graphs with the uniform distribution as the initial distribution.  We
have $\P[X_t=X_0]\,-\,1/|V_o|$ is
\be
&\P[X_t=X_0]\,-\,\frac{1}{|V_o|}\;\;=\;\;\frac{1}{|V_o|}\suli{j=2}{|V_o|}
e^{-t(1-\beta_j)} \;\;\ge\;\;\frac{1}{|V_o|} e^{-t\lambda} \;\;\ge
\;\;\frac{1}{|V_o|} e^{-\frac{t\,K}{|V_o|^{1/2}}}.&\;\;\hfill
\qed \ee

{\em Theorem \ref{theo:3}; Lower bound}: Compare this with
[\ref{mathieu}], Lemma 2.2 and [\ref{woess}].  Let $G$ be transitive,
with a unimodular, transitive subgroup of $\textrm{Aut}(G)$, the
automorphism group of $G$. Given $\w\in \Omega$, for $G'(\w)$ being
the whole percolation subgraph of $G$, the graph $G_o$ is the
connected subgraph of $G'(\w)$ induced by $\Cn(\w)$,
i.e. $V_o=\Cn(\w)$. (In what follows, we will drop the dependence on
$\w$, wherever it doesn't cause confusion. For example, we write $\Cn$
instead of $\Cn(\w)$.)\\

From Theorem \ref{theo:upperbound}, iii.), since $G_o$ is almost
surely finite, there is a lower bound for the expected return
probability of the delayed random walk. Namely, since
\be
\E_\mu\left[\P[X_t=X_0]\right]\;\;-\;\;\E_\mu\left[\frac{1}{|\Cn|}
  \right]\;\;\ge\;\;\E_\mu\left[\frac{1}{|\Cn|}\,e^{-\frac{t
      K}{\sqrt{|\Cn|}}}\chi_{{}_{|\Cn|>288}}\right],
\ee
and due to the assumption $t > \sqrt{288}$, we have 
\be
\E_\mu\left[\frac{1}{|\Cn|}\,e^{-\frac{t K}{\sqrt{|\Cn|}}}\chi_{{}_{|\Cn|\ge t^2}}
\right]\;\;\ge\;\;\suli{m\ge
  t^2}{\infty}\frac{1}{m}\,e^{-\frac{tK}{\sqrt{m}}}\phi(m)
\;\;\ge\;\;e^{-K}\suli{m\ge t^2}{\infty} \frac{1}{m} \,\phi(m).  
\ee

The lower bound of Theorem \ref{theo:3} now follows by Lemma
\ref{lemma:heavytail}, with $D=e^{-K}\frac{A/2}{1+(2B/A)^{1/b}}$ and
by applying Lemma \ref{lemma:mtp} to express $P_t(o)$ by the
normalised trace.  \hfill \qed\\

{\em Theorem \ref{theo:3}; Upper bound}: By assumption, $\mu$ is
invariant under a unimodular transitive subgroup of Aut$(G)$, and by
the remark after Corollary \ref{cor:z2_trees} there are almost surely
only finite cluster. In particular $\mu$-a.s.  $|\Cn|<\infty$. \\

Let $N=|\Cn|$. In order to use Theorem \ref{theo:upperbound} most
effectively, we want to choose $k\in\{1, ..., N-2\}$ as small as
possible while keeping the exponents of the same order in $t/N^2$. We
differentiate between two cases: First, we assume $\lfloor
N\sqrt{q\lambda}\rfloor+1 \le N-2$, where $q=\d^2(\d+2)/32$. Then, we
choose $k$ in Theorem \ref{theo:upperbound}, i.) such that
\be 
\lambda\;:=\;1\;-\;\beta_2\; \le\; \frac{32}{\d^2(\d+2)} \cdot
\frac{k^2}{N^2}.
\ee

This is accomplished if we set $k=\lfloor N\sqrt{q\lambda}\rfloor+1$.
(Note, $k\le N-2$.) This choice implies $N\sqrt{q \lambda}< k \le
1+N\sqrt{q \lambda}$.  Setting $c=\sqrt{\pi q/2}$, it follows
\bel 
\P[X_t \;=\;X_0]\;\;&\le&\;\;
\frac{1}{N}\;+\;\left((\frac{2}{N}\,+\,2\sqrt{q\lambda}) \;+\;
\frac{c}{\sqrt{t}}\right)e^{-\lambda t}.  \label{eq:secondterm}
\eel

From $e^{-x}\le y^y/x^y$ and $e^{-x}\le ((y-1/2)/x)^{y-1/2}$ for $y>
\eh$, we get
\be 
\P[X_t \;=\;X_0]\;\;&\le&\;\;
  \frac{1}{N}\;+\;\frac{1}{t^y}\left(y^y\left(\frac{2}{N\l^y}
  +\frac{2\sqrt{q}}{\lambda^{y-1/2}}\right) \;+\; c
  \frac{(y-\eh)^{y-1/2}}{\lambda^{y-1/2}}\right).  
\ee

Now using the Poincar\'e inequality $\lambda\ge \d/(4 N^2)$, we obtain
the following estimate:
\bel
\P[X_t =X_0]&\le&
\frac{1}{N}+\frac{1}{t^y}\left(y^y\left(2\frac{16^y N^{2y-1}}{\d^y}
+\frac{\sqrt{q}\,2^{2y} N^{2y-1}}{\d^{y-1/2}}\right) + c
\frac{(y-\eh)^{y-1/2}(2N)^{2y-1})}{\d^{y-1/2}}\right)\nn\\
&\le & \frac{1}{N}\;+\;c_\d\frac{N^{2y-1}}{t^y},\label{eq:Nexp}
\eel

with\\[-0.6cm]
\bel
c_\d\;\;=\;\;2^{2y}\left(\frac{y}{\d}\right)^y\left( 2^{2y+1}\;+\;
\frac{\d\sqrt{\d(\d+2)}}{4\sqrt{2}}\left(1\,+\, 
\sqrt{2\pi} \right)\right).\label{eq:constant}
\eel

We have used that $b>0$ and that by (\ref{eq:assump1}) and the remark
after Corollary \ref{cor:z2_trees} the exponent of $N$ in
(\ref{eq:Nexp}) is $\a:=2y-1<b$, so that $1/2<y<1/2+b/2$ and
$(y-\eh)^{y-\eh}< 2 y^{\,y}$. With $\sqrt{\d(\d+2)}\le \d+2$ and
$(1+\sqrt{2\pi})/(4 \sqrt{2})\le 1/2$, this leads to the upper bound
$c_\d \le \left(4\a/\d\right)^\a(2\cdot 4^\a \,+\, \d(\d+2)/2)$.
Since $b\le 2$, we see that the constant in front of $N^{2y-1}/t^y$ in
(\ref{eq:Nexp}) is bounded below by $\eh$, independently of $\d$.

Now, turning to $\lfloor N\sqrt{q\lambda}\rfloor+1 \ge N-2$, which is
equivalent to $\lambda\ge (1- 3/N)^2/q\ge 1/(16 q)$, if $N\ge 4$.  If
$N<4$, we have the Poincar\'e inequality $\lambda\ge \d/(4N^2)\ge
\d/36$.  So, in both cases, the function $t\mapsto\P[X_t=X_0]-1/N$ is
decreasing exponentially fast. Since it is smaller than $1$, the
overall estimate covering all three possibilities (including the
polynomially decreasing one) is given if the constant $c_\d > 1$ in
(\ref{eq:constant}) is multiplied by five, yielding $5\cdot
\left(4b/\d\right)^b(2\cdot 4^b \,+\, \d(\d+2)/2)$. 
\\[-0.3cm]

Taking the expected value of both sides of the inequality and applying
Lemma \ref{lemma:mtp} to express $P_t(o)$ by the normalised trace
yields the result.  \hfill \qed\\

{\em Corollary \ref{cor:z2_trees}, i.); Upper bound}: Since Bernoulli
bond percolation on the Euclidean lattice is invariant under the
unimodular transitive group of translations of the Euclidean lattice,
this is a special case of Theorem \ref{theo:3}. The result follows
from the well-known fact [\ref{kesten2}], that there exists $\a>0$,
s.t. $\E_\mu[|\Cn|^\a]<\infty$.  \hfill \qed\\

{\em Corollary \ref{cor:z2_trees}, i.); Lower bound}: By the power law
inequality $\Phi(m)=\P_{\mu}[|\Cn| \ge m]\ge \eh\, m^{-\eh}$ (see
[\ref{grimmett}], Theorem 11.89), we have $a=\eh$ in
(\ref{eq:assump1}).  For any $t_o>0$ it is now possible to choose
$C_2$ depending on $t_o$ such that $C_2^{-1}(288)^{-(1+1/\a)}=1$.  So,
by choosing $t_o=1$ the given estimate follows from Theorem
\ref{theo:3} for $\a < b$ and for all $t\ge 1$. \hfill \qed\\

{\em Corollary \ref{cor:z2_trees}, ii.); Upper and Lower bound}: It is
well-known that for the homogeneous tree of finite degree, $a=b=\eh$
(see [\ref{aldous}], [\ref{kolchin}], and [\ref{timar}]). Just as in
the previous proof, the constant $C_\d$ can be chosen so large, that
the estimate is valid for all $t\ge 1$. \hfill \qed\\[-0.3cm]

{\em  Theorem \ref{theo:lower}}:
\hspace{0.7cm}a.) Let
$I_t=\E_\mu\left[\P[X_t=X_0]-\frac{1}{|\Cn|}\right]$, and let
$\l=1-\beta_2$ be the smallest non-zero eigenvalue of $\id-P$, as
above. We have, for any $c>0$, $I_t\ge
\E_\mu\left[\;\P[X_t=X_0]-\frac{1}{|\Cn|}\; \right|\left.\; \l \le
  \;\frac{c}{|\Cn|^2}\right]\cdot\P\left[\;\l \le
  \;c/|\Cn|^{2}\;\right]$, and that by (\ref{eq:trace}) \be
\E_\mu\left[\;\P[X_t=X_0]-\frac{1}{|\Cn|}\; \right|\left.\; \l
  \le\;\frac{c}{|\Cn|^2}\right] \;\ge\;
\E_\mu\left[\;\frac{e^{-t\l}}{|\Cn|}\;\;\right|\left.\; \l
  \le\;\frac{c}{|\Cn|^2}\right]\;=\; \E_\mu\left[\;
  \frac{e^{-\frac{c t}{|\Cn|^2}}}{|\Cn|}\right]\,.  \ee

b.) Let for $\w\in\Omega$ the diameter $D(\w)$ of $\Cn(\w)$ be
defined by $ D\;\;=\;\;\max_{v, w \in \Cn} d(v, w)$,

with $d(., .)$ the graph metric of $G_o$. Let $\pi=(v_0, v_1, v_2,
... , v_D)$ be a geodesic path in $G_o$ of length $D$. Consider the
function $g:\Cn\to \R$ with $g(v) = \cos(\pi k/D)$ where $k$ is
uniquely defined by 
$d(v_k, v)=\min\{\,d(v_j,v)\;|\; j\in\{0, ..., D\}\,\}$. \\

Now, we show that if for some number $\e>0$ it holds $\e |\Cn|\le D$,
then the function $g$ gives an upper estimate of
$\l$ in terms of $|\Cn|^{-2}$:
\be
\lambda \;=\;\min\limits_{f\perp \textrm{const}} 
\frac{\sum_{i< j\in\Cn}(f_i - f_j)^2}{\sum_{v\in\Cn}|f(v)|^2}\;\le\;
\frac{\sum_{v\sim w\in\Cn}(g(v) - g(w))^2}{\sum_{v\in\Cn}|g(v)|^2}\;\le\;
\frac{\sum_{j=1}^D(g(v_j) - g(v_{j-1}))^2}{\sum_{j=1}^D|g(v_j)|^2},
\ee

where the second inequality results from neglecting the terms in the
denominator not belonging to the geodesic $\pi$. By Taylor's Theorem $\cos(\pi j/D)
= \cos(\pi (j-1)/D) + (\pi/D)\sin(\pi (j-1)/D) + O(1/D^2)$ as
$D\mapsto \infty$, so  for some number $c>0$
\be
\lambda \;\;\le\;\; 
\frac{\pi^2}{D^2}\frac{\sum_{j=1}^D(\sin(\pi (j-1)/D)^2}{\sum_{j=1}^D|
\cos(\pi j/D)|^2}\left(1 + O(\frac{1}{D^2})\right)\;\;\le\;\;
\frac{c}{D^2}\;\;\le\;\;\frac{c}{\e^2|\Cn|^2}.
\ee

c.) By Markov's inequality, for $\a<b$
\be
\P_\mu\left[\frac{|\Cn|}{D}\ge \e^{-1}\right]\;\;\le\;\;
\e^{\a}\E_\mu\left[\frac{|\Cn|^\a}{D^\a}\right]\;\;\le\;\;
\e^{\a}\E_\mu[|\Cn|^\a]
\ee

of which the right-hand side can be made smaller than one by choosing
$\e$ sufficiently small. For such an $\e$ the probability of the
complement is positive, or, in other words, $C:=\P_\mu[\e|\Cn| < D]>
0$. So, from b.),  $\P[\l<c/(\e^2|\Cn|^2)]$ for some $c>0$
with a probability bounded below by $C>0$.\\

d.) Let $\phi(m)=\P_\mu[|\Cn|=m]$, and $t>0$.  Under the assumptions
\be 
\sum_{m>\sqrt{t}}\,\frac{\phi(m)}{m}\;\ge\; A
\sum_{m>\sqrt{t}}\,m^{-a-2} \; \ge\; A \int_{\sqrt{t}}^\infty
x^{-a-2}dx \;=\; \frac{A}{a+1}\frac{1}{(\sqrt{t})^{a+1}} 
\ee
and so, by the foregoing arguments (a., b., c.), $I_t$ is bounded from below by
\be
C\cdot\E_\mu\left[\frac{e^{-t/(\e^2|\Cn|^2)}}{|\Cn|}\right]\;\ge\;
C \suli{m>\sqrt{t}}{}\frac{1}{m}e^{-\frac{t}{\e^2 m^2}}\P[|\Cn|=m]\;\ge\;
\frac{C\, A \,e^{-1/\e^{2}}}{(a+1)} \,\,t^{-\frac{a+1}{2}}.\;\;\;\qed\\[-0.3cm]
\ee

{\em Corollary \ref{cor:three_fourths}}:\hspace{0.3cm} 
Since by Corollary \ref{cor:z2_trees}, ii.) it holds for all $\e>0$
that
\be
\lim\limits_{t\to\infty}
\frac{\log(\;P_t(o)\;-\;\E_\mu[1/|\Cn|]\;)}{\log t}
\;\;\;\le\;\;\;-\frac{3}{4}\;\;+\;\;\e, \ee 

it must be true for $\e=0$, and the upper bound follows.  Furthermore,
it is well known [\ref{fies}] that for critical percolation on the
homogeneous tree $\P_\mu[|\Cn|=m] \sim m^{-3/2}$.  Therefore, the
assumptions of Theorem \ref{theo:lower} are fulfilled where $a=b=1/2$
(see [\ref{grimmett}], Chap. 10.1, and [\ref{hughes}], Chap. 1.3), which
implies the lower bound. \hfill\qed\\

{\em Theorem \ref{theo:2}; Upper bound}: The integrated density of states
$N(E)$ obeys [\ref{kirsch}, \ref{muesto}, \ref{sob}] the relation
$\int_0^\infty e^{-tE}dN(E)=\E_\mu[\P_o[X_t=o]]$, such that by Theorem
\ref{theo:3}, i.)\\[-0.6cm]
\be e^{-t\e}(N(\e)\;-\;N(0))\;\;\le\;\; \int_0^\e
e^{-tE}dN(E)\;\;\le\;\; P_t\;-\;\kappa
\;\;\le\;\;c_{4}\,\E_\mu[|\Cn|^\a]\,t^{-\nu},\\[-0.5cm]\nn \ee
where $\nu=\eh(1+\a)$, with $\a$ such that $\E_\mu[|\Cn|^\a]<\infty$, and
$c_4=\left(8+\sqrt{3\pi}\right)$. Choosing
$t=\nu/\e$ and thereby optimising the upper bound for $N(\e)\;-\;N(0)$
leads to the result.\\

{\em Theorem \ref{theo:2}; Lower bound}: Again, by $\int_0^\infty
e^{-tE}dN(E)=\E_\mu[\P_o[X_t=o]]$, Lemma \ref{lemma:mtp} and Corollary 
\ref{cor:z2_trees}, with $\a>0$ s.t. $\E_\mu[|\Cn|^\a]<\infty$, 
\be 
\frac{C_2^{-1}}{t^{(1+1/\a)}}\le \int_0^\infty e^{-t\,E}dN(E)
\;\;\le\;\; \int_0^\e dN(E) \;+\;e^{-t\e}\int_\e^\infty dN(E)
\;\;\le\;\;N(\e)-N(0)+e^{-t\e}.  \\[-0.3cm]\nn
\ee

So, $N(\e)-N(0)\ge \eh C_2^{-1}t^{-(1+1/\a)} -
e^{-t\e}$. Choosing $t=-(\bar{c}/\e)\log \e$ for $\e>0$ produces the
result if, for example,
$\bar{c}=2\cdot (1+1/\a)$. Then 
$C_3=\max\{ C_2^{-1}/(\bar{c}\log\e)^{1+\a^{-1}},$
$ c_4 \;\E_\mu[|\Cn|a]\}$. \hfill \qed\\

{\em Theorem \ref{theo:num}}:\; Bernoulli bond percolation on the
$d$-dimensional Euclidean lattice is a percolation invariant under the
unimodular translation group of the lattice.  The degree is $\d=2\cdot
d$.  Assuming subcritical Bernoulli bond-percolation, we have
existence of the first moment of the cluster size.  By repeating the
argument of the proof of the upper bound in Theorem \ref{theo:3}
(which lead to (\ref{eq:secondterm})) with Theorem
\ref{theo:upperbound} ii.)  instead of \ref{theo:upperbound} i.)
yields for all $t>0$ with $q=4/(d^2(d+1))$
\be
\E_\mu\P[X_t \;=\;X_0]\;\;&\le&\;\;\E_\mu
\left[\frac{1}{|\Cn|}\right]\;+\; \E_\mu\left[\frac{2
    k}{|\Cn|}e^{-\frac{t}{ |\Cn|^2}}\;\;+\;\; \frac{2}{q\,t}
  \frac{|\Cn|}{k} \exp\left(-\frac{q\,t\, k^2}{|\Cn|^2} \right)
  \right].
\ee

Now, choosing $k=1$, and using $\exp x\;\le\; 1/x$ for $x>0$ gives
for all $t>0$
\be
\E_\mu\P[X_t \;=\;X_0]\;\;&\le&\;\;\E_\mu
\left[\frac{1}{|\Cn|}\right]\;+\; \E_\mu\left[2 \frac{ |\Cn|}{t}
\;\;+\;\; \frac{2 |\Cn|}{q\,t}  \right] .  
\ee

Calling $\kappa(p)=\E_\mu[1/|\Cn|]$ (note the difference to
[\ref{grim_num}] regarding the cluster which consist of only one
vertex), letting $\chi(p):=\E_\mu[|\Cn|]$ and noting $2 +
2/q=(d^3+d^2+4)/2$, leads to the lower bound after a subsequent
application of Lemma \ref{lemma:mtp}, and a rearrangement of the terms
in the inequality. \\

The other constant $\frac{20}{d}(4+d(d+1))$ follows from the method
used for proving the upper bound of Theorem \ref{theo:3}, and by using
$b=1$ and setting $\a$ in $\E_p[|\Cn|^\a]<\infty$ equal to $b$, which is
possible due to the existence of the first moment.\\

The upper bound follows from the observation $P_t\;-\; \kappa(p) =$
$\E_\mu[(1/|\Cn|)\cdot \Sp \exp(-t(1-P))] \ge 0 $, since $1-P$ has
only non-negative eigenvalues. \hfill \qed

\section{Acknowledgment} I am grateful for the discussion with Peter
M\"uller concerning the asymptotics of $P_t$ and its comparison with
the incipient infinite cluster in the case of homogeneous
trees. Thanks also to Adam Timar informing me about [\ref{timar}], to
Geoffrey Grimmett and an `anonymous expert' for remarks leading to the
discussion following Corollary \ref{cor:three_fourths}, and to Daniel
Lenz and Brian Rider for their helpful comments regarding improvements
of the style of the presentation.  This work has been written with the
support of the Project P18703 of the Austrian Science Foundation
(FWF), and during the author's stay at the Mathematisches Institut of
the University of Jena.

\section{Bibliography}

\begin{enumerate}

\item{D. Aldous: `Asymptotic fringe disttributions for general
    families of random trees', Ann. of Probab. Vol1. 1, No.2,
    228-266, (1991)}\label{aldous}

\item{S. Alexander, R. Orbach: `Density of states on fractals:
  ``fractons'' ', J. Physique (Paris) Lett. 43, 625-631,
  (1982)}\label{alob}

\item{T. Antunovic, I. Veselic: `Spectral asymptotics of percolation
hamiltonians on amenable Cayley graphs', In: Proceedings of OTAMP
2006. Operator Th.: Adv. Appl. 2007}\label{Tonci}

\item{A. Bandyopadhyay, J. Steif, A. Timar: `On the Cluster Size
  Distribution for Percolation on Some General Graphs',
  arXiv:0805.3620v1, Theorem 2.c, (2008)}\label{timar}

\item{M. T. Barlow, T. Kumagai: `Random walk on the incipient infinite
  cluster on trees', Illinois Journ. Math., 50, 33-65, (2006)
}\label{baku}

\item{V. Batagelj, T. Pisanski: `Hamiltonian cycles in the cartesian
  product of a tree and a cycle', Discr. Math,. 38, 311-312,
  (1982)}\label{bapi}

\item{I. Benjamini, R. Lyons; Y. Peres, O. Schramm: `Critical
  percolation on any nonamenable group has no infinite cluster',
  Ann. Probab.  27, no. 3, 1347-1356, (1999)}\label{blps2}

\item{A.G. Boshier: `Enlarging properties of graphs', Ph.D. thesis,
  Royal Holloway and Bedford New College, University of London,
  (1987)}\label{boshier}

\item{V. V. Dimakopoulos, L. Palios, A. S. Poulakidas: `On the
  hamiltonicity of the cartesian product' Inf. Process. Lett. 96(2):
  49-53, (2005)}\label{dipapo}

\item{M. E. Fisher, J. W. Essam: `Some cluster size and percolation
  problems', Jour. Math. Phys. 2, 609-627, (1961)}\label{fies}

\item{L.R.Font\`es, P. Mathieu: `On symmetric random walks with random
  rates. Probability Theory and Related Fields', vol. 134, No 4,
  565-602, (2006)}\label{mathieu}

\item{G.R.Grimmett: `Percolation', 2nd edition, Springer
  (1999)}\label{grimmett}

\item{G.R.Grimmett: `On the number of clusters in the percolation
  model', J.London Math.Soc. (2), 13 (1976), 346-350}\label{grim_num}

\item{J. van den Heuvel: `Hamiltonian Cycles and Eigenvalues of
  Graphs', Lin. Alg. Appl. 226-228:723-730, (1995)}\label{heuvel}

\item{B. D. Hughes: `Random Walks and Random Environments', Volume 2,
Clarendon Press, Oxford 1995}\label{hughes}

\item{H. Kesten: `The Critical Probability of Bond Percolation on the
  Square Lattice Equals 1/2', Commun. Math. Phys. 74, 41-59 (1980)
}\label{kesten3}

\item{H. Kesten: `Subdiffusive behaviour of random walk on a random
  cluster', Annales de l'Institut Henri Poincar\'e , Probab. et
  Statistiques 22, 425-487, (1986)}\label{kesten}

\item{H. Kesten: `Scaling relations for 2D-percolation',
    Comm. Math. Phys. 109, 09-156, (1987)}\label{kesten2}

\item{W.Kirsch, P.M\"uller: `Spectral properties of the Laplacian on
  bond-percolative graphs', Math. Z. 252, 899-916 (2006) 
}\label{kirsch}

\item{V.F. Kolchin: `Random Mappings', Optimization Software
    Inc. Publ. Div., New York, (Transl. of Russion Original),
    (1986)}\label{kolchin}

\item{G. Kozma, A. Nachmias: `The Alexander-Orbach conjecture holds in
  high dimensions', Invent. math. (2009) 178: 635-654}\label{kona}

\item{R. Lyons, Y. Peres: `Probability on trees and networks',
  web-book: 

    http://mypage.iu.edu/$\sim$rdlyons/prbtree/prbtree.html
}\label{lyons}

\item{B. Mohar:`Isoperimetric numbers of graphs', Jour. Comb. Theory,
  B 47, 274-291, (1989)}\label{mohar}

\item{P. M\"uller, P. Stollmann: `Spectral asymptotics of the
  Laplacian on supercritical bond-percolation graphs', 
  Jour. Func. Anal. 252, Is. 1, 1 11, 233-246, (2007)}\label{muesto}

\item{J.R. Norris: `Markov Chains', Cambridge University Press,
  (1997), Chap. 2}\label{norris}

\item{L.Saloff-Coste: `Lectures on finite Markov Chains', Saint-Flour
  summer school: lecture notes, LNM 1665, (1997)}\label{laurent}

\item{F. Sobieczky: `An interlacing technique for spectra of random
    walks and its application to finite percolation cluster',
    J. Theoret. Probab. 23 (2010), no. 3, 639-670}\label{sob}

\item{W. Woess: `Random Walks on Infinite Graphs and Groups',
  Cambridge Tracts in Mathematics 138, Camb. Univ. Press, (2000),
  (Chap. 12.13)}\label{woess}

\end{enumerate}

\end{document}